\documentclass[12pt]{article}%
\usepackage{amsfonts}
\usepackage{sw20bams}
\usepackage{amsmath}
\usepackage{amssymb}
\usepackage{graphicx}%
\providecommand{\U}[1]{\protect\rule{.1in}{.1in}}
\begin{document}

\title{Covariance within Random Integer Compositions}
\author{Steven Finch}
\date{October 19, 2020}
\maketitle

\begin{abstract}
Fix a positive integer $N$. \ Select an additive composition $\xi$ of $N$
uniformly out of $2^{N-1}$ possibilities. \ The interplay between the number
of parts in $\xi$ and the maximum part in $\xi$ is our focus. \ It is not
surprising that correlations $\rho(N)$ between these quantities are negative;
we earlier gave inconclusive evidence that $\lim_{N\rightarrow\infty}\rho(N)$
is strictly less than zero. \ A proof of this result would imply asymptotic
dependence.\ \ We now retract our presumption in such an unforeseen outcome.
\ Similar experimental findings apply when $\xi$ is a $1$-free composition,
i.e., possessing only parts $\geq2.$\ 

\end{abstract}

\footnotetext{Copyright \copyright \ 2020 by Steven R. Finch. All rights
reserved.}An \textbf{unrestricted} additive composition of $N$ is a sequence
of positive integers, called parts, that sum to $N$. \ The number of
unrestricted compositions is $2^{N-1}$. \ For example, the compositions of $5$
are%
\[%
\begin{array}
[c]{cccc}%
\{5\}, & \{4,1\}, & \{3,2\}, & \{2,3\},\\
\{1,4\}, & \{3,1,1\}, & \{2,2,1\}, & \{2,1,2\},\\
\{1,3,1\}, & \{1,2,2\}, & \{1,1,3\}, & \{2,1,1,1\}\\
\{1,2,1,1\}, & \{1,1,2,1\}, & \{1,1,1,2\}, & \{1,1,1,1,1\}.
\end{array}
\]
If a composition of $N$ is chosen uniformly at random from all possibilities,
then \cite{SF-tcs5, F1-tcs5, F2-tcs5}
\[
1+m_{n}=\mathbb{E}(\text{number of parts})=1+\dfrac{n}{2},
\]%
\[
s_{n}^{2}=\mathbb{V}(\text{number of parts})=\dfrac{n}{4}%
\]
where $n=N-1$, $\mathbb{E}$ denotes mean and $\mathbb{V}$ denotes variance
(uncorrected for bias);%
\[
1+\mu_{n}=\mathbb{E}(\text{maximum part})=1+\frac{1}{2^{n}}\left[
z^{n}\right]
{\displaystyle\sum\limits_{k=1}^{\infty}}
\left(  \frac{1}{1-2z}-\frac{1-z^{k}}{1-2z+z^{k+1}}\right)  ,
\]%
\[
\sigma_{n}^{2}=\mathbb{V}(\text{maximum part})=\frac{1}{2^{n}}\left[
z^{n}\right]
{\displaystyle\sum\limits_{k=1}^{\infty}}
(2k-1)\left(  \frac{1}{1-2z}-\frac{1-z^{k}}{1-2z+z^{k+1}}\right)  -\mu_{n}^{2}%
\]
where $[z^{n}]$ denotes the coefficient of $z^{n}$ in the subsequent Taylor
series expansion. \ For example, if $N=5$, then $m_{4}=2$, $s_{4}^{2}=1$,
$\mu_{4}=27/16$ and $\sigma_{4}^{2}=247/256$. \ Up to small periodic
fluctuations \cite{Byd-tcs5, Sch-tcs5, Alx-tcs5}, we have%
\[
\mathbb{E}(\text{maximum part})\sim\dfrac{\ln(n)}{\ln(2)}+\left(
\dfrac{\gamma}{\ln(2)}-\dfrac{1}{2}\right)  ,
\]%
\[
\mathbb{V}(\text{maximum part})\sim\dfrac{1}{12}+\dfrac{\pi^{2}}{6\ln(2)^{2}}%
\]
asymptotically as $N\rightarrow\infty$.

A \textbf{restricted} additive composition of $N$ obeys an extra condition
that no parts are equal to $1$. \ A\ more descriptive name is \textbf{1-free}.
\ The number of restricted compositions is $d_{n}$, the $n^{\text{th}}$
Fibonacci number, where $d_{0}=0$, $d_{1}=1$ and $d_{n}=d_{n-1}+d_{n-2}$.
\ For example, the $1$-free compositions of $7$ are%
\[%
\begin{array}
[c]{cccc}%
\{7\}, & \{5,2\}, & \{4,3\}, & \{3,4\},\\
\{2,5\}, & \{3,2,2\}, & \{2,3,2\}, & \{2,2,3\}.
\end{array}
\]
If a $1$-free composition of $N$ is chosen uniformly at random from all
possibilities, then \cite{F1-tcs5, F2-tcs5}%
\[
1+m_{n}=\mathbb{E}(\text{number of parts})=1+\frac{1}{d_{n}}\left[
z^{n}\right]  \frac{z^{3}}{\left(  1-z-z^{2}\right)  ^{2}},
\]%
\[
s_{n}^{2}=\mathbb{V}(\text{number of parts})=\frac{1}{d_{n}}\left[
z^{n}\right]  \frac{z^{3}\left(  1-z+z^{2}\right)  }{\left(  1-z-z^{2}\right)
^{3}}-m_{n}^{2},
\]%
\[
1+\mu_{n}=\mathbb{E}(\text{maximum part})=1+\frac{1}{d_{n}}\left[
z^{n}\right]
{\displaystyle\sum\limits_{k=1}^{\infty}}
\left(  \frac{1-z^{2}}{1-z-z^{2}}-\frac{1-z^{2}-z^{k}+z^{k+1}}{1-z-z^{2}%
+z^{k+1}}\right)  ,
\]%
\[
\sigma_{n}^{2}=\mathbb{V}(\text{maximum part})=\frac{1}{d_{n}}\left[
z^{n}\right]
{\displaystyle\sum\limits_{k=1}^{\infty}}
(2k-1)\left(  \frac{1-z^{2}}{1-z-z^{2}}-\frac{1-z^{2}-z^{k}+z^{k+1}}%
{1-z-z^{2}+z^{k+1}}\right)  -\mu_{n}^{2}.
\]
For example, if $N=7$, then $m_{6}=5/4$, $s_{6}^{2}=7/16$, $\mu_{6}=13/4$ and
$\sigma_{6}^{2}=27/16$. \ Up to small periodic fluctuations, we conjecture
that \cite{F2-tcs5}%
\[
\mathbb{E}(\text{maximum part})\sim\dfrac{\ln(n)}{\ln(\varphi)}+\left(
\dfrac{\gamma}{\ln(\varphi)}-1\right)  ,
\]%
\[
\mathbb{V}(\text{maximum part})\sim\dfrac{1}{12}+\dfrac{\pi^{2}}{6\ln
(\varphi)^{2}}%
\]
as $N\rightarrow\infty$, where $\varphi=(1+\sqrt{5})/2$ is the Golden mean.

What is missing among these results? \ We have not yet conveyed any sense of
how the two quantities are interrelated. \ Define a correlation coefficient%
\begin{align*}
\rho(N)  &  =\frac{\mathbb{E}\left[  (\text{number of parts})(\text{maximum
part})\right]  -(1+m_{n})(1+\mu_{n})}{s_{n}\sigma_{n}}\\
&  =\frac{\mathbb{E}\left[  (\text{number of }1\text{s})(\text{longest run of
}0\text{s})\right]  -m_{n}\mu_{n}}{s_{n}\sigma_{n}}%
\end{align*}
where the first expression (involving integer compositions) is understandable
but the second expression (involving bitstrings) needs clarification.
\ Avoiding details for now, let us simply provide some numerical data in Table 1.

\begin{center}%
\[%
\begin{tabular}
[c]{|c|c|c|}\hline
$n$ & $\rho(N)$ for unrestricted case & $\rho(N)$ for $1$-free case\\\hline
$100$ & $-0.441772$ & $-0.530911$\\\hline
$200$ & $-0.361888$ & $-0.439875$\\\hline
$300$ & $-0.319761$ & $-0.391011$\\\hline
$400$ & $-0.292051$ & $-0.358533$\\\hline
$500$ & $-0.271797$ & $-0.334641$\\\hline
$600$ & $-0.256049$ & $-0.315973$\\\hline
$700$ & $-0.243295$ & $-0.300791$\\\hline
$800$ & $-0.232656$ & $-0.288084$\\\hline
$900$ & $-0.223581$ & $-0.277216$\\\hline
$1000$ & $-0.215704$ & $-0.267762$\\\hline
$1100$ & $-0.208773$ & $-0.259428$\\\hline
$1200$ & $-0.202606$ & $-0.251998$\\\hline
$1300$ & $-0.197066$ & $-0.245313$\\\hline
$1400$ & $-0.192050$ & $-0.239249$\\\hline
\end{tabular}
\ \ \ \ \
\]
Table 1:\ Correlation between $($number of parts$)$ and $($maximum part$)$
within random integer compositions as a function of $N=n+1$.\bigskip
\end{center}

\noindent Acceleration of convergence is possible for each sequence,
suggesting (without proof \cite{F2-tcs5}) that limits are nonzero as
$N\rightarrow\infty$. \ We now must retract such presumptive and unjustified
thoughts. \ A more careful study leads to a revised conjecture:%
\[%
\begin{array}
[c]{ccc}%
\rho(N)\sim C\ln(N)^{-5/2} &  & \text{for some }C<0
\end{array}
\]
thus in particular $\rho(N)\rightarrow0^{-}$, consistent with asymptotic
independence. \ This is plainly what intuition leads everyone to foresee. \ A
rigorous treatment would be good to see someday.

\section{Unconstrained and Pinned Solus Bitstrings}

Given a random unconstrained bitstring of length $n=N-1$, we have%
\[%
\begin{array}
[c]{ccc}%
\mathbb{E}(\text{number of }1\text{s})=n/2, &  & \mathbb{V}(\text{number of
}1\text{s})=n/4
\end{array}
\]
because a sum of $n$ independent Bernoulli($1/2$) variables is Binomial($n$%
,$1/2$). \ Expressed differently, the average density of $1$s in a string is
$1/2$, with a corresponding variance $1/4$. \ The word \textquotedblleft
unconstrained\textquotedblright\ offers that, in the sampling process, all
$2^{n}$ strings are included and equally weighted.

If we append the string with a $1$, calling this $\eta$, then there is a
natural way \cite{HS-tcs5} to associate $\eta$ with an additive composition
$\xi$ of $N$. \ For example, if $N=10$,%
\[
\eta=0110100111\longleftrightarrow\xi=\{2,1,2,3,1,1\}
\]
i.e., parts of $\xi$ correspond to \textquotedblleft waiting
times\textquotedblright\ for each $1$ in $\eta$. \ The number of parts in
$\xi$ is equal to the number of $1$s in $\eta$ and the maximum part in $\xi$
is equal to the duration of the longest run of $0$s in $\eta$, plus one.

In this paper, the word \textquotedblleft constrained\textquotedblright%
\ refers to the logical conjunction of two requirements:

\begin{itemize}
\item A bitstring is \textbf{pinned} if its first bit is $0$ and its last bit
is $0$.

\item A bitstring is \textbf{solus} if all of its $1$s are isolated.
\end{itemize}

\noindent The latter was discussed in \cite{F1-tcs5, F2-tcs5}; additionally
imposing the former is new. \ Given a random pinned solus bitstring of length
$n=N-1$, formulas for $\mathbb{E}($number of $1$s$)$ and $\mathbb{V}($number
of $1$s$)$ are best expressed using generating functions.

If we append the string with $1$ to construct $\eta$, then the associated
$\xi$ is a composition of $N$ with all parts $\geq2$. \ For example, if
$N=15$,%
\[
\eta=010001010010101\longleftrightarrow\xi=\{2,4,2,3,2,2\}.
\]
It should now be clear why, starting with the original $n$-bitstring,%
\[
1+\mathbb{E}(\text{number of }1\text{s})=\mathbb{E}(\text{number of parts}),
\]%
\[
1+\mathbb{E}(\text{longest run of }0\text{s})=\mathbb{E}(\text{maximum part})
\]
for both scenarios, but the corresponding variances are always equal.

Nej \& Satyanarayana Reddy \cite{NS-tcs5} gave an impressive recursion for the
number $F_{n}(x,y)$ of unconstrained bitstrings of length $n$ containing
exactly $x$ $0$s and a longest run of exactly $y$ $0$s:
\[
F_{n}(x,y)=\left\{
\begin{array}
[c]{lll}%
{\displaystyle\sum\limits_{i=\kappa}^{y-1}}
F_{n-i-1}(x-i,y)+%
{\displaystyle\sum\limits_{j=0}^{y}}
F_{n-y-1}(x-y,j) &  & \text{if }1\leq x\leq n-2\text{ and }\varepsilon
_{n}(x,y)=1,\\
\lambda_{n}(y) &  & \text{if }x=n-1\text{ and }\varepsilon_{n}(x,y)=1,\\
0 &  & \text{otherwise,}%
\end{array}
\right.
\]%
\[%
\begin{array}
[c]{ccc}%
F_{n}(0,0)=1-\kappa, &  & F_{n}(n,n)=1
\end{array}
\]
where $n\geq x\geq y$ (of course) and $\kappa=0$,
\[
\varepsilon_{n}(x,y)=\left\{
\begin{array}
[c]{lll}%
1 &  & \text{if }n\geq x\text{ and }\left\lfloor \dfrac{n}{n-x+1}\right\rfloor
\leq y\leq x,\\
0 &  & \text{otherwise}%
\end{array}
\right.
\]
and
\[
\lambda_{n}(y)=\left\{
\begin{array}
[c]{lll}%
1 &  & \text{if }n\text{ is odd and }y=\dfrac{n-1}{2},\\
2 &  & \text{otherwise.}%
\end{array}
\right.
\]
Consequently, the numerator of $\mathbb{E}\left[  (\text{number of }%
1\text{s})(\text{longest run of }0\text{s})\right]  $ for $n$-bitstrings is
\[
\left\{
{\displaystyle\sum\limits_{x=0}^{n}}
\,%
{\displaystyle\sum\limits_{y=0}^{x}}
(n-x)y\,F_{n}(x,y)\right\}  _{n=1}^{\infty}%
=\{0,2,11,40,122,338,881,2202,5337,12634,\ldots\};
\]
equivalently, the numerator of $\mathbb{E}\left[  (\text{number of
parts})(\text{maximum part})\right]  $ for $N$-compositions is%
\[
\left\{
{\displaystyle\sum\limits_{x=0}^{N}}
\,%
{\displaystyle\sum\limits_{y=0}^{x}}
(N-x)(y+1)F_{n}(x,y)\right\}  _{N=1}^{\infty}%
=\{1,4,14,42,115,296,732,1757,4125,9516,\ldots\}.
\]
The denominator is $2^{n}$. \ Returning to the unrestricted example, the
covariance for $N=5$ is $\tfrac{40}{16}-(2)\left(  \tfrac{27}{16}\right)
=\tfrac{115}{16}-(1+2)\left(  1+\tfrac{27}{16}\right)  $. \ Correlations for
selected small $N$ turn out to be%
\[
\left\{  \rho(N):N=5,11,21,51\right\}
=\{-0.890799,-0.752444,-0.654958,-0.530128\}
\]
and Table 1 exhibits values for larger $N=101,201,\ldots$.\newpage

By a similar argument, we deduce the number $G_{n}(x,y)$ of pinned solus
bitstrings of length $n$ containing exactly $x$ $0$s and a longest run of
exactly $y$ $0$s. \ The recursion is identical to before (with $F$ replaced by
$G$). \ The initial conditions appear alike, but here we have $\kappa=1$.
\ Also, a different $\varepsilon_{n}(x,y)$ applies:
\[
\varepsilon_{n}(x,y)=\left\{
\begin{array}
[c]{lll}%
1 &  & \text{if }n\geq x\text{ and }\left(  \left\lfloor \dfrac{n}%
{n-x+1}\right\rfloor \leq y<x\text{ or }x=y=n\right)  ,\\
0 &  & \text{otherwise}%
\end{array}
\right.
\]
and a different $\lambda_{n}(y)$:%
\[
\lambda_{n}(y)=\left\{
\begin{array}
[c]{lll}%
1 &  & \text{if }n\text{ is odd and }y=\dfrac{n-1}{2},\\
2 &  & \text{if }\left\lfloor \dfrac{n-1}{2}\right\rfloor <y<n-1,\\
0 &  & \text{otherwise.}%
\end{array}
\right.
\]
Consequently, the numerator of $\mathbb{E}\left[  (\text{number of }%
1\text{s})(\text{longest run of }0\text{s})\right]  $ under constraints is
\[
\left\{
{\displaystyle\sum\limits_{x=0}^{n}}
\,%
{\displaystyle\sum\limits_{y=0}^{x}}
(n-x)y\,G_{n}(x,y)\right\}  _{n=1}^{\infty}%
=\{0,0,1,4,10,26,54,118,230,458,864,1632,\ldots\};
\]
equivalently, the numerator of $\mathbb{E}\left[  (\text{number of
parts})(\text{maximum part})\right]  $ under restrictions is%
\[
\left\{
{\displaystyle\sum\limits_{x=0}^{N}}
\,%
{\displaystyle\sum\limits_{y=0}^{x}}
(N-x)(y+1)G_{n}(x,y)\right\}  _{N=1}^{\infty}%
=\{0,2,3,8,17,34,70,131,255,466,868,1565,\ldots\}.
\]
The denominator is $d_{n}$. \ Returning to the $1$-free example, the
covariance for $N=7$ is $\tfrac{26}{8}-\left(  \tfrac{5}{4}\right)  \left(
\tfrac{13}{4}\right)  =\tfrac{70}{8}-\left(  1+\tfrac{5}{4}\right)  \left(
1+\tfrac{13}{4}\right)  $. \ Correlations for selected small $N$ turn out to
be%
\[
\left\{  \rho(N):N=7,11,21,51\right\}
=\{-0.945611,-0.860467,-0.763395,-0.629068\},
\]
i.e., dependency is more significant than earlier. \ Table 1 exhibits values
for larger $N=101,201,\ldots$.

\section{Sketches of Proofs I}

Let $\Omega$ be a set of finite bitstrings and $\Omega_{n}^{x,y}$ be the
subset of $\Omega$ consisting of strings of length $n$ containing exactly $x$
$0$s and a longest run of exactly $y$ $0$s. \ Let $\Omega_{n,0}^{x,y}$ and
$\Omega_{n,1}^{x,y}$ be the subset of $\Omega_{n}^{x,y}$ of strings starting
with $0$ and $1$ respectively. \ 

Assume that $\Omega$ consists of all unconstrained strings. \ If $\omega
\in\Omega_{n,1}^{x,y}$, then $\omega$ is of the form $1\omega_{1}$ where
$\omega_{1}\in\Omega_{n-1,0}^{x,y}\cup\Omega_{n-1,1}^{x,y}$. \ If $\omega
\in\Omega_{n,0}^{x,y}$, then $\omega$ is either of the form%
\[%
\begin{array}
[c]{ccccccc}%
\underset{i}{\underbrace{00\ldots0}}\omega_{2} &  & \text{where} &  &
\omega_{2}\in\Omega_{n-i,1}^{x-i,y} & \text{and} & 1\leq i\leq y-1
\end{array}
\]
or%
\[%
\begin{array}
[c]{ccccccc}%
\underset{y}{\underbrace{00\ldots0}}\omega_{3} &  & \text{where} &  &
\omega_{3}\in\Omega_{n-y,1}^{x-y,j} & \text{and} & 0\leq j\leq y.
\end{array}
\]
We have%
\[
\left\vert \Omega_{n,1}^{x,y}\right\vert =\left\vert \Omega_{n-1,0}%
^{x,y}\right\vert +\left\vert \Omega_{n-1,1}^{x,y}\right\vert =\left\vert
\Omega_{n-1}^{x,y}\right\vert ,
\]%
\begin{equation}
\left\vert \Omega_{n,0}^{x,y}\right\vert =%
{\displaystyle\sum\limits_{i=1}^{y-1}}
\left\vert \Omega_{n-i,1}^{x-i,y}\right\vert +%
{\displaystyle\sum\limits_{j=0}^{y}}
\left\vert \Omega_{n-y,1}^{x-y,j}\right\vert =%
{\displaystyle\sum\limits_{i=1}^{y-1}}
\left\vert \Omega_{n-i-1}^{x-i,y}\right\vert +%
{\displaystyle\sum\limits_{j=0}^{y}}
\left\vert \Omega_{n-y-1}^{x-y,j}\right\vert
\end{equation}
hence%
\[
\left\vert \Omega_{n}^{x,y}\right\vert =%
{\displaystyle\sum\limits_{i=0}^{y-1}}
\left\vert \Omega_{n-i-1}^{x-i,y}\right\vert +%
{\displaystyle\sum\limits_{j=0}^{y}}
\left\vert \Omega_{n-y-1}^{x-y,j}\right\vert
\]
upon addition. \ This proof of the recurrence for $F_{n}(x,y)$ appeared in
\cite{NS-tcs5}. \ 

Assume instead that $\Omega$ consists of all solus strings. \ If $\omega
\in\Omega_{n,1}^{x,y}$, then $\omega$ is of the form $1\omega_{1}$ where
$\omega_{1}\in\Omega_{n-1,0}^{x,y}$. \ We have%
\[
\left\vert \Omega_{n,1}^{x,y}\right\vert =\left\vert \Omega_{n-1,0}%
^{x,y}\right\vert =\left\vert \Omega_{n-1}^{x,y}\right\vert -\left\vert
\Omega_{n-1,1}^{x,y}\right\vert ,
\]
that is,%
\[
\left\vert \Omega_{n}^{x,y}\right\vert =\left\vert \Omega_{n,1}^{x,y}%
\right\vert +\left\vert \Omega_{n+1,1}^{x,y}\right\vert =\left\vert
\Omega_{n-1,0}^{x,y}\right\vert +\left\vert \Omega_{n,0}^{x,y}\right\vert .
\]
From formula (1) in the preceding,%
\[
\left\vert \Omega_{n,0}^{x,y}\right\vert =%
{\displaystyle\sum\limits_{i=1}^{y-1}}
\left\vert \Omega_{n-i,1}^{x-i,y}\right\vert +%
{\displaystyle\sum\limits_{j=0}^{y}}
\left\vert \Omega_{n-y,1}^{x-y,j}\right\vert =%
{\displaystyle\sum\limits_{i=1}^{y-1}}
\left\vert \Omega_{n-i-1,0}^{x-i,y}\right\vert +%
{\displaystyle\sum\limits_{j=0}^{y}}
\left\vert \Omega_{n-y-1,0}^{x-y,j}\right\vert
\]
which gives a recurrence underlying what we called $\tilde{F}_{n}(x,y)$ in
\cite{F2-tcs5}.

Let us turn attention to various boundary conditions.\ \ For either
unconstrained or solus strings,%
\[
\Omega_{n}^{n-1,n-1}=\left\{  1\underset{n-1}{\underbrace{00\ldots0}%
},\;\underset{n-1}{\underbrace{00\ldots0}}1\right\}  ;
\]
if $n$ is odd, then%
\[
\Omega_{n}^{n-1,(n-1)/2}=\left\{  \underset{(n-1)/2}{\underbrace{00\ldots0}%
}1\underset{(n-1)/2}{\underbrace{00\ldots0}}\right\}  ;
\]
if $n$ is even, then%
\[
\Omega_{n}^{n-1,n/2}=\left\{  \underset{(n-2)/2}{\underbrace{00\ldots0}%
}1\underset{n/2}{\underbrace{00\ldots0}},\;\underset{n/2}{\underbrace{00\ldots
0}}1\underset{(n-2)/2}{\underbrace{00\ldots0}}\right\}  .
\]
These imply the expression for $\lambda_{n}(y)$. For pinned strings, the
latter two results hold, but the former becomes $\Omega_{n}^{n-1,n-1}%
=\varnothing$. \ The expression for $\varepsilon_{n}(x,y)$ comes from
\cite{NS-tcs5}:%
\[
F_{n}(x,y)>0\Longleftrightarrow\left\{
\begin{array}
[c]{lll}%
x+\left\lfloor \dfrac{x}{y}\right\rfloor \leq n &  & \text{if }y>0\text{ and
}y\nmid x\text{,}\\
x+\dfrac{x}{y}-1\leq n &  & \text{if }y>0\text{ and }y\mid x;
\end{array}
\right.
\]
\
\[
G_{n}(x,y)>0\Longleftrightarrow\left\{
\begin{array}
[c]{lll}%
x+\left\lfloor \dfrac{x}{y}\right\rfloor \leq n &  & \text{if }y>0\text{ and
}y\nmid x\text{,}\\
x+\dfrac{x}{y}-1\leq n &  & \text{if (}x>y>0\text{ or }x=y=n\text{) and }y\mid
x.
\end{array}
\right.
\]
For completeness' sake, we give the analog of Table 1 for pinned and solus strings.

\begin{center}%
\[%
\begin{tabular}
[c]{|c|c|c|}\hline
$n$ & $\rho$ for pinned case & $\rho$ for solus case\\\hline
$100$ & $-0.445112$ & $-0.525562$\\\hline
$200$ & $-0.363340$ & $-0.437637$\\\hline
$300$ & $-0.320638$ & $-0.389680$\\\hline
$400$ & $-0.292661$ & $-0.357617$\\\hline
$500$ & $-0.272255$ & $-0.333956$\\\hline
$600$ & $-0.256411$ & $-0.315434$\\\hline
$700$ & $-0.243592$ & $-0.300351$\\\hline
$800$ & $-0.232906$ & $-0.287715$\\\hline
$900$ & $-0.223795$ & $-0.276900$\\\hline
$1000$ & $-0.215891$ & $-0.267488$\\\hline
$1100$ & $-0.208938$ & $-0.259187$\\\hline
$1200$ & $-0.202753$ & $-0.251783$\\\hline
$1300$ & $-0.197198$ & $-0.245119$\\\hline
$1400$ & $-0.192170$ & $-0.239074$\\\hline
\end{tabular}
\ \ \ \ \ \ \ \
\]
Table 2:\ Correlation between $($number of $1$s$)$ and $($longest run of
$0$s$)$ within random bitstrings as a function of $n$.\bigskip
\end{center}

\section{Sketches of Proofs II}

Let $\Omega$ and $\Omega_{n}^{x,y}$ be as before. \ Assume that $\Omega$
consists of all pinned strings. \ If $\omega\in\Omega_{n}^{x,y}$, then
$\omega$ is either of the form
\begin{equation}%
\begin{array}
[c]{ccc}%
\underset{\geq1}{\underbrace{00\ldots0}}\underset{\geq2}{\underbrace{11\ldots
1}}0\ldots, &  & i=0
\end{array}
\end{equation}
or%
\begin{equation}%
\begin{array}
[c]{ccc}%
\underset{i}{\underbrace{00\ldots0}}10\ldots, & \;\;\;\; & 1\leq i\leq y.
\end{array}
\end{equation}
On the one hand, the subset of strings $\omega$ satisfying (2)\ corresponds to
$\Omega_{n-1}^{x,y}$ upon deleting the leftmost $1$. \ On the other hand, the
subset of strings $\omega$ satisfying (3)\ corresponds to
\[%
\begin{array}
[c]{ccccc}%
\Omega_{n-i-1}^{x-i,y} &  & \text{if} &  & i<y
\end{array}
\]
and to
\[%
\begin{array}
[c]{ccccccccc}%
\Omega_{n-y-1}^{x-y,j} &  & \text{for some} &  & 0\leq j\leq y &  & \text{if}
&  & i=y
\end{array}
\]
upon deleting the leftmost block of $0$s and the leftmost $1$. \ These
observations give%
\[
\left\vert \Omega_{n}^{x,y}\right\vert =%
{\displaystyle\sum\limits_{i=0}^{y-1}}
\left\vert \Omega_{n-i-1}^{x-i,y}\right\vert +%
{\displaystyle\sum\limits_{j=0}^{y}}
\left\vert \Omega_{n-y-1}^{x-y,j}\right\vert
\]
which interestingly is the same recurrence as that for $F_{n}(x,y)$. \ Clearly
however%
\[%
\begin{array}
[c]{ccc}%
\left\vert \Omega_{n}^{0,0}\right\vert =0, &  & \left\vert \Omega_{n}%
^{n,n}\right\vert =1
\end{array}
\]
are the initial conditions here.

Assume instead that $\Omega$ consists of all pinned solus strings. \ The
condition described by (2)\ is no longer satisfied by any $\omega\in\Omega
_{n}^{x,y}$ because $1$s are now isolated. \ Hence the case $i=0$ is removed
from the summation, implying that $\kappa=1$, which gives the recurrence for
$G_{n}(x,y)$.

\section{Acknowledgements}

R, Mathematica and Maple have been useful throughout. I am grateful to Ernst
Joachim Weniger, Claude Brezinski and Jan Mangaldan for very helpful
discussions about convergence acceleration. \ Dr. Weniger's software code and
numerical computations were especially appreciated. \

\end{document}